\input amstex
\documentstyle{amsppt}
\topmatter
\magnification=\magstep1
\pagewidth{5.2 in}
\pageheight{6.7 in}
\abovedisplayskip=10pt
\belowdisplayskip=10pt
\NoBlackBoxes
\title
$q$-Euler numbers and polynomials associated with basic zeta
functions
\endtitle
\author  Taekyun Kim  \endauthor
\affil\rm{{Division of General Education-Mathematics,}\\
{ Kwangwoon University, Seoul 139-701, S. Korea}\\
{e-mail: tkkim$\@$kw.ac.kr}}
\endaffil

\define\C{\Bbb C_p}

\define\Z{\Bbb Z_p}

\define\ra{\rightarrow}

\abstract{ We consider the $q$-analogue of Euler zeta function
which is defined  by
$$\zeta_{q, E}(s)=[2]_q\sum_{n=1}^{\infty}\frac{(-1)^nq^{ns}}{[n]_q^s}, \text{ $0<
q<1$, $\Re(s)> 1$.}$$ In this paper, we give the $q$-extension of
Euler  numbers which  can be viewed as interpolating of the above
$q$-analogue of Euler zeta function at negative integers , in the
same way that Riemann zeta function interpolates Bernoulli numbers
at negative integers. Also, we will treat some identities of the
$q$-extension of the Euler numbers by using fermionic $p$-adic
$q$-integration on $\Bbb Z_p $.}
\endabstract
\thanks 2000 Mathematics Subject Classification  11S80, 11B68 \endthanks
\thanks Key words and phrases: $q$-Bernoulli number, $p$-adic $q$-integrals, zeta function, Dirichlet series \endthanks
\rightheadtext{   } \leftheadtext{ $q$-analogue of Riemann zeta
function  }
\endtopmatter

\document

\head 1. Introduction \endhead
Throughout this paper $\Bbb Z_p$, $\Bbb Q_p$, $\Bbb C$ and $\Bbb C_p$ will respectively denote the ring of $p$-adic rational integers,
the field of $p$-adic rational numbers, the complex number field and the completion of algebraic closure of $\Bbb Q_p $.

The $p$-adic absolute value in $\Bbb C_p$ is normalized so that
$|p|_p=\frac1p .$ When one talks of $q$-extension, $q$ is
variously considered as an indeterminate, a complex number $q\in
\Bbb C$ or a $p$-adic number $q\in \Bbb C_p.$ If $q\in \Bbb C ,$
then we normally assume $|q|< 1 ,$ and when $q\in\Bbb C_p,$ then
we normally assume $|q-1|_p < 1$.  We use the notation :
$$[x]_q=[x:q]= \frac{1-q^x}{1-q}, \text{  and  } [x]_{-q}=\frac{1-(-q)^x}{1+q}.$$
Note that $\lim_{q \rightarrow 1} [x]_q = x$ for $ x \in \Bbb Z_p$
in presented $p$-adic case.

Let $UD(\Bbb Z_p )$ be denoted by the set of uniformly differentiable functions on $\Bbb Z_p $.

For $ f \in UD( \Bbb Z_p) $, let us start with the expression
$$ \frac{1}{[p^N ]_{-q}} \sum_{ 0 \leq j < p^N } (-q)^j f(j) = \sum_{0 \leq j < p^N} f(j)\mu_{-q} ( j + p^N \Bbb Z_p )$$
representing analogue of Riemann's sums for $f$, cf.[1-30].

The fermionic $p$-adic $q$-integral of $f$ on $\Bbb Z_p$ will be
defined as the limit $( N \rightarrow \infty )$ of these sums,
which it exists. The fermionic $p$-adic $q$-integral of a function
$ f \in UD(\Bbb Z_p)$ is defined as
$$ \int_{\Bbb Z_p} f(x) d \mu_{-q}(x) = \lim_{N \rightarrow \infty} \frac{1}{[p^N]_{-q}}
\sum_{0 \leq j < p^N } f(j) (-q)^j,\text{ (see [5, 6, 16] )}.$$ For
$d$ a fixed positive integer with $(p,d)=1$, let
$$X=X_d=\varprojlim_N \Bbb Z/dp^N\Bbb Z , \;\;X_1=\Bbb Z_p,$$
$$X^*=\bigcup\Sb 0<a<dp\\ (a,p)=1\endSb a+dp\Bbb Z_p,$$
$$a+dp^N\Bbb Z_p=\{x\in X\mid x\equiv a\pmod{dp^N}\},$$
where $a\in \Bbb Z$ lies in $0\leq a<dp^N$, (see [1-30]).

Let $\Bbb N$ be the set of positive integers.
 For $m, k\in N$,
 the $q$-Euler polynomials
$E_{m}^{(-m,k)}(x, q)$ of higher order in the variables $x$
 in $\C$ by making use of the $p$-adic $q$-integral
, cf.[5, 6], are defined by
$$\multline
E_{m,q}^{(-m,k)}(x)=\undersetbrace\text{$k$
times}\to{\int_{\Z}\int_{\Z}\cdots\int_{\Z}}
[x+x_1 +x_2 +\cdots +x_k ]_q^m \\
\cdot
q^{-x_1(m+1)-x_2(m+2)-\cdots-x_k(m+k)}d\mu_{-q}(x_1)d\mu_{-q}(x_2)\cdots
d\mu_{-q}(x_k).
\endmultline \tag1$$

Now, we define the $q$- Euler numbers of higher order as follows:
$$E_{m,q}^{(-m, k)} =E_{m,q}^{(-m, k)} (0).$$

From (1), we can derive
$$\aligned
&E_{m, q}^{(-m,k)}\\
&=\lim_{N \ra
\infty}\frac{1}{[p^N]_{-q}^k}\sum_{x_1=0}^{p^N-1}\cdots
\sum_{x_k=0}^{p^{N}-1}
[x_1+\cdots+x_k]_q^m(-1)^{x_1+\cdots+x_k} q^{-x_{1} m-\cdots-x_k(m+k-1)}\\
&=\frac{[2]_q^k}{(1-q)^m}\sum_{i=0}^{m} \binom mi (-1)^i
\frac{1}{(1+q^{i-m})(1+q^{i-m-1})\cdots (1+q^{i-m-k+1})},
\endaligned \tag 2$$
where $\binom{m}{i}$ is binomial coefficient.

 Note that
$\lim_{q\rightarrow 1}E_{m,q}^{(-m,k)}=E_m^{(k)}$ where
$E_m^{(k)}$ are ordinary Euler numbers of order $k$, which are
defined as
$$\left(\frac{2}{e^t+1}\right)^k=\sum_{n=0}^{\infty}E_n^{(k)}\frac{t^n}{n!}.$$

 By (1), (2), it is easy to see that
$$E_{m,q}^{(-m, 1)}(x)=\sum_{i=0}^{m}\binom mi q^{xi}E_{i,q}^{(-m, 1)} [x]^{m-i}
=\frac{ [2]_q}{(1-q)^m}\sum_{j=0}^m q^{jx}\binom mj (-1)^j
\frac{1}{1+q^{j-m}}. \tag 3 $$

We  define the $ q$-analogue of Euler zeta function which is
defined  as
$$ \zeta_{q, E}(s)=[2]_q \sum_{n=1}^{\infty}
\frac{(-1)^nq^{sn}}{[n]_q^s}, \text{ $ q\in\Bbb R$ with $0<q<1$
and $s\in\Bbb C$. }\tag 4$$ The numerator ensures the convergence.
In (4), we can consider the following problem:

`` Are there $q$-Euler numbers which  can be viewed as
interpolating of $\zeta_{q, E}(s)$ at negative integers, in the
same way that Riemann zeta function interpolates Bernoulli numbers
at negative integers "?

In this paper, we give the value $\zeta_{q, E}(-m)$, for $
m\in\Bbb N$, which is a answer of the above problem and construct
a new complex $q$-analogue of Hurwitz's type Euler zeta function
and $q$-$L$-series related to $q$-Euler numbers. Also, we will
treat some interesting identities of $q$-Euler numbers.

\head 2. some identities of  $q$-Euler numbers $E_{m,q}^{(-m,
1)}$.
\endhead

In this section, we assume $q\in \Bbb C_p$ with $|1-q|_p < 1.$ By
(1), we see that
$$\aligned
E_{n,q}^{(-n, 1)}(x)&= \int_{ X}q^{-(n+1)t}[x+t]_q^n d
\mu_{-q}(t)\\
&=\frac{[2]_q}{[2]_{q^d}}[d]_q^{n}\sum_{i=0}^{d-1}(-1)^iq^{-ni}\int_{\Bbb
Z_p}
 q^{-(n+1)dt}[\frac{x+i}{d}+t]_{q^d}^n d\mu_{-q^d }(t) .\endaligned$$
Thus we have
$$ E_{n,q}^{(-n,1)}(x)=\frac{[2]_q}{[2]_{q^d}}[d]_q^{n}\sum_{i=0}^{d-1} (-1)^iq^{-ni}E_{n,q^d}^{(-n, 1)}
( \frac{ x+i}{d} ), \tag5$$ where $d$, $n$ are positive integers
with $d \equiv 1 (\mod 2)$.

If we take $ x=0 $, then we have
$$E_{m,q}^{(-m,1)}=\frac{[2]_q}{[2]_{q^n}}\sum_{k=0}^m
\binom{m}{k}[n]_q^kE_{k,q^n}^{(-m,1)}\sum_{j=0}^{n-1}(-1)^jq^{-(m-k)j}[j]_q^{m-k},
\text{ where $n\equiv 1(\mod 2)$} . \tag6$$

 From (6), we can easily derive  the following equation (7).
 $$E_{m,q}^{(-m, 1)}-\frac{[2]_q}{[2]_{q^n}}[n]_q^m E_{m, q^n}^{(-m, 1)}
=\frac{[2]_q}{[2]_{q^n}}\sum_{k=0}^{m-1} \binom mk [n]_q^k
E_{k,q^n}^{(-m, 1)} \sum_{j=1}^{n-1}q^{-(m-k)j} (-1)^j[j]_q^{m-k} .
\tag7$$ It is easy easy to see that $\lim_{q \rightarrow 1}
E_{m,q}^{(-m, 1)} =E_m$, where $E_m $ are the $m$-th ordinary Euler
numbers, cf. [5]. From (7),  we note that
$$ ( 1-n^m ) E_m =
\sum_{k=0}^{m-1} \binom mk n^k E_k \sum_{j=1}^{n-1}(-1)^j
j^{m-k}.$$

Let $F_q(t,x)$ be generating function of $E_{n,q}^{(-n, 1)}$ as
follows:
$$F_q(t,x)=\sum_{k=0}^{\infty}E_{k,q}^{(-k,1)}(x) \frac{t^k}{k!}. \tag 8$$
By (3), (8), we easily see:
$$\aligned
 F_q(t,x)&=[2]_q\sum_{k=0}^{\infty}\left(\frac{1}{(q-1)^k}\sum_{j=0}^k(-1)^{k-j}\binom{k}{j}\frac{q^{jx}}{1+q^{j-k}}
 \right)\frac{t^k}{k!}\\
 &=\sum_{k=0}^{\infty}\left([2]_q\sum_{n=0}^{\infty}(-1)^nq^{-kn}[n+x]_q^{k}\right)
\frac{t^k }{k!}. \endaligned\tag 9$$

Differentiating both sides with respect to $t$ in (5), (6) and
comparing coefficients, we obtain the following: \proclaim{
Theorem 1} For $m \geq 0,$  we have
$$ E_{m, q}^{(-m, 1)}(x)=[2]_q\sum_{n=0}^{\infty} q^{-nm}[n+x]_q^{m}(-1)^n. \tag10$$
\endproclaim

\proclaim{ Corollary 2} Let $m \in \Bbb N.$  Then there exists
$$ E_{m, q}^{(-m, 1)}=[2]_q\sum_{n=1}^{\infty} q^{-nm}[n]_q^{m}(-1)^n,
\text{ and  } E_{0,q}^{(0, 1)}= \frac{[2]_q}{2}. \tag11$$
\endproclaim
 Note that Corollary 2 is a $q$-analogue of $ \zeta_{E}( m )$,
for any positive integer $m$.

Let $\chi$ be a primitive Dirichlet character with conductor $ d
\in \Bbb N $ with $d\equiv 1 (\mod 2)$.

For $ m \in \Bbb N$, we define
$$ E_{ m, \chi, q}^{(-m, 1)}= \int_{X}
q^{-(m+1)x}\chi(x)[x]^m d \mu_{-q}(x), \text{ for }  m \geq
0.\tag12
$$

Note that
$$\aligned
 E_{m, \chi,q}^{(-m,
 1)}&=\frac{[2]_q}{[2]_{q^d}}[d]_q^m\sum_{i=0}^{d-1}q^{-mi}(-1)^i\chi(i)\int_{\Bbb
 Z_p}q^{-dx(m+1)}[\frac{i}{d}+x]_{q^d}^m d\mu_{-q^d}(x)\\
  &=\frac{[2]_q}{[2]_{q^d}}[d]_q^{m}\sum_{i=0}^{d-1}
\chi(i)(-1)^iq^{-mi}E_{m, q^d}^{(-m, 1)}(\frac {i}{d})
 .\endaligned\tag13$$

\head 3. $q$-analogs of zeta functions \endhead In this section,
we assume $q\in \Bbb R$ with $0<q<1 .$ Now we consider the
$q$-extension of the Euler zeta function as follows:
$$ \zeta_{q, E}(s)=[2]_q\sum_{n=1}^{\infty}(-1)^n\frac{q^{ns}}{[n]_q^s}, \text{ where $s \in \Bbb C$ .}$$
By(11), we obtain the following theorem.

\proclaim{ Theorem 3} For $m \in \Bbb N, $  we have
$$ \zeta_{q, E}(-m)=E_{m,q}^{(-m, 1)}.$$
\endproclaim

From Theorem 1, we can also define the $q$-extension of Hurwitz's
type Euler $\zeta$-function as follows: For $s \in \Bbb C ,$
define
$$ \zeta_{q, E}(s, x)=[2]_q\sum_{n=0}^{\infty}\frac{(-1)^nq^{sn}}{[n+x]_q^s}. \tag 14$$
Note that $\zeta_{q, E}(s,x)$ is an analytic continuation in whole
complex $s$-plane .

By (14) and Theorem 1, we have the following theorem.

 \proclaim{Theorem 4} For any positive integer $k$, we have
$$ \zeta_{q, E}(-k, x) =E_{k,q}^{(-k, 1)}(x,q). \tag 15$$
\endproclaim

For $d\in \Bbb N$ with $d\equiv 1 (\mod 2)$, let $\chi$ be
Dirichlet character with conductor $ d $. By (13), the generalized
q-Euler numbers attached to $\chi$ can be defined as
$$ E_{m, \chi,q}^{(-m, 1)} =\frac{[2]_q}{[2]_{q^d}}[d]_q^{m}\sum_{i=0}^{d-1}
\chi(i)q^{-mi}(-1)^i E_{m, q^d}^{(-m, 1)}(\frac {i}{d}) .\tag16 $$

For $ s \in \Bbb C$, we define
$$L_{q, E} (s, \chi)=[2]_q\sum_{n=1}^{\infty}\frac{ \chi(n)(-1)^n q^{sn}}{[n]_q^s}. \tag 17$$
It is easy to see that
$$L_{q, E}(\chi, s)=\frac{[2]_q}{[2]_{q^d}} [d]_q^{-s}\sum_{a=1}^d\chi(a)(-1)^aq^{sa}\zeta_{q^d, E}(s, \frac ad). \tag18$$
By (16), (17), (18), we obtain the following theorem.

 \proclaim{
Theorem 5}Let $k \in \Bbb N$. Then there exists
$$ L_{q, E}( -k, \chi) =E_{k, \chi, q}^{(-k, 1)}.$$
\endproclaim
Let $a$ and $F$ be integers with $ 0< a < F $. For $ s \in \Bbb C
$, we consider the functions $H_q(s, a, F)$ as follows:
$$ H_{q, E}( s,a,F)=[2]_q\sum_{m\equiv a ( F), m>0} \frac{q^{ms}(-1)^m}{[m]_q^s}
= \frac{[2]_q}{[2]_{q^F}}[F]_q^{-s}(-1)^a q^a\zeta_{q^F}(s, \frac
aF).$$ Then we have
$$ H_{q, E}( -n, a, F )=(- 1)^aq^a\frac{ [2]_q}{[2]_{q^F}}[F]_q^{n}E_{n, q^F}^{(-n, 1)}(\frac aF ),$$
  where $n$ is any positive integer .

In the recent paper, the $q$-analogue of Riemann zeta function
related to twisted $q$-Bernoulli numbers was studied by Y. Simsek
(see [1, 26, 30]). In [30], Y. Simsek have studies  the twisted
$q$-Bernoulli numbers which can be viewed as an interpolating of
 the $q$-analogue of Riemann zeta function  at negative integers.
 In this paper, we have shown that the $q$-analogue of Euler zeta function
interpolates the $q$-Euler numbers at negative integers, in the same
way that Riemann zeta function interpolates Bernoulli numbers at
negative integers, cf. [5, 6, 30].

\enddocument